\documentclass[11pt,a4paper]{amsart}
\usepackage{amsmath}

\usepackage{amssymb,latexsym}

\newcommand{\dmath}{\displaystyle}

\newtheorem{theorem}{Theorem}
\newtheorem{corollary}{Corollary}
\newtheorem{lem}{Lemma}
\newenvironment{dem}{\textbf{Proof}}{}
\begin{document}

\title[minimality of $x/\|x\|$]{An induction principle for the weighted $p$-energy minimality of $x/\|x\|$}

\author{ Jean-Christophe Bourgoin }

\address{ Laboratoire de Math\'ematiques et Physique Th\'eorique,
UMR CNRS 6083, Universit\'e de Tours, Parc de Grandmont, F-37200
Tours France}

\email{}

\begin{abstract}

In this paper, we investigate minimizing properties of the map $x/\|x\|$ from the Euclidean unit ball $\mathbf{B}^{n}$ to its boundary $\mathbb{S}^{n-1}$, for the weighted energy functionals $E^n_{p,\alpha}(u)=\int_{\mathbf{B}^{n}}  \|x\|^{\alpha}\|\nabla u\|^p dx$.  We establish the following induction principle: if the map $\frac{x}{\|x\|}:\mathbf{B}^{n+1}\to \mathbb{S}^n$ minimizes $E^{n+1}_{p,\alpha}$ among the maps $u: \mathbf{B}^{n+1}\to \mathbb{S}^n$ satisfying $u(x)=x$ on
$\mathbb{S}^n$, then the map $\frac{y}{\|y\|}:\mathbf{B}^n\to\mathbb{S}^{n-1}$ minimizes $E^{n}_{p,\alpha+1}$ among the maps $v: \mathbf{B}^n\to\mathbb{S}^{n-1}$ satisfying $v(y)=y$ on $\mathbb{S}^{n-1}$. 

This result enables us to enlarge the range of values of $p$ and $\alpha$ for which $x/\|x\|$ minimizes $E^n_{p,\alpha}$.  
\end{abstract}

\maketitle

\section {Introduction and statement of main results}\label{1}
Let $\mathbf{B}^n$ be the unit Euclidean ball of $\mathbb{R}^n$ and $\mathbb{S}^{n-1}$ its bourdary.
For any couple $ (\alpha,p)$ of real numbers, with $\alpha\ge0$ and $p\ge 1$,
we define the $r^\alpha$-weighted $p$-energy functional of a map $u~: \mathbf{B}^n\rightarrow \mathbb{S}^{n -1}$ by
\[
E^n_{p,\alpha}(u)=\int_{\mathbf{B}^n}\|x\|^{\alpha}\|\nabla u\|^p dx.
\]
This functional is nothing but the $p$-energy functional associated with the Riemannian metric $r^{\alpha(n-p)}g_{euc}$ on the ball $\mathbf{B}^n$, where $g_{euc}$ is the Euclidean metric, the sphere $\mathbb{S}^{n-1}$ being endowed with its standard metric. The functional $E^n_{p,\alpha}$ is to be considered on the Sobolev space
\[ W_{\alpha}^{1,p}(\mathbf{B}^{n},\mathbb{S}^{n-1})=\{u \in W^{1,p}((\mathbf{B}^n r^{\alpha(n-p)}g_{euc}),\mathbb{R}^n) ; \|u\|=1 a.e\}.
\]
The question is to know which map minimizes $E_{p,\alpha}$ among the maps in $W_{\alpha}^{1,p}(\mathbf{B}^{n},\mathbb{S}^{n-1})$ satisfying $u(x)=x$ on $\mathbb{S}^{n-1}$.

Following the arguments of Hildebrant, Kaul and Widman \cite{HKW}, the map $x/\|x\|$ is a natural candidate to be the minimizer of $E^n_{p,\alpha}$, for $p\in [1,n+\alpha)$ (notice that $x/\|x\|\in W_{\alpha}^{1,p}(\mathbf{B}^{n},\mathbb{S}^{n-1})$ if and only if $p<n+\alpha$).

The minimality of $x/\|x\|$ was first established for the 2-energy functional $E_2^n:=E_{2,0}^n$ by J\"ager and Kaul \cite{JK} in dimension $n \geq 7$, then by Brezis, Coron and Lieb \cite{BCL} in dimension 3. Coron and Gulliver \cite{CG} actually proved the minimality of $x/\|x\|$ for the $p$-energy functional $E_p^n:=E_{p,0}^n$ for any integer $p \in \{1, \cdots, n-1\}$ and any dimension $n \geq 3$. An alternative proof using the null Lagrangian method (or calibration method) of this last result was obtained by Lin \cite{L} and Avellaneda and Lin \cite{AL}.

Hardt and Lin  \cite{HL1, HL2} and Hardt, Lin and Wang \cite{HLW1, HLW2} developed the study of the  singularities of $p$-harmonic and $p$-minimizing maps and obtained results extending those of Schoen and Uhlenbelck \cite{SU1, SU2}. A consequence of their results is the minimality of $x/\|x\|$ for $E_p^n$ for any $p \in\left(n-1,n\right)$. Finally, Hong \cite{Ho1} and Wang \cite{W} have proved independently the minimality of $x/\|x\|$ for $E_p^n$ in dimension $n \geq 7$ for any $p \leq n-2\sqrt{n-1}$.

In order to close the question of the $E^n_p$-minimality of  $x/\|x\|$ for the remaining values of $p$, Hong \cite{Ho2} suggested a new idea. Indeed, he observed that, for any $p\in \left(2,n\right)$, the minimality of $x/\|x\|$ for $E^n_p$ follows from the minimality of $x/\|x\|$ for the $r^{2-p}$-weighted 2-energy $E^n_{2,2-p}$. Unfortunately, we have proved in a previous paper \cite{B1}, the existence of an interval of values of $p\in \left(1,n-1\right)$ for which the map $x/\|x\|$ fails to be a minimizer of $E^n_{2,2-p}$.
Nevertheless, we proved that $x/\|x\|$ minimizes the $r^\alpha$-weighted $p$-energy $E^n_{p,\alpha}$ for any $\alpha\ge0$ and any integer $p\le n-1$. Actually, we obtained in \cite{B2} the minimality of $x/\|x\|$ for more general weighted $p$-energy functionals of the form $E^n_{p,f}(u)=\int_{\mathbf{B}^n}f(\|x\|)\|\nabla u\|^p dx$, where $p\le n-1$ is an integer and $f: [0,1] \rightarrow \mathbb{R}^+ $ is a continuous non-decreasing function.

Our aim in this paper is to show that the minimizing properties of  $x/\|x\|$ in dimension $n$ and $n+1$ are not independent. Indeed, we will prove that, for any $p \in [1,n+\alpha+1)$, if $x/\|x\|$ minimizes the $r^\alpha$-weighted $p$-energy $E^{n+1}_{p,\alpha}$ in dimension $n+1$, then it also minimizes the $r^{\alpha+1}$-weighted $p$-energy $E^{n}_{p,\alpha+1}$ in dimension $n$.

 \begin{theorem}\label{th1}

Let $n \geq 2$ be an integer, $\alpha \ge 0$ be a real number and $p \in [1,n+\alpha+1)$. If the map $\frac{x}{\|x\|}:\mathbf{B}^{n+1}\mapsto \mathbb{S}^n$ minimizes the  $r^\alpha$-weighted $p$-energy $E^{n+1}_{p,\alpha}$ among the maps $u \in W_{\alpha}^{1,p}
(\mathbf{B}^{n+1}, \mathbb{S}^n)$ satisfying $u(x)=x$ on
$\mathbb{S}^n$, then the map $\frac{y}{\|y\|}:\mathbf{B}^n\mapsto\mathbb{S}^{n-1}$ minimizes the $r^{\alpha+1}$-weighted $p$-energy $E^{n}_{p,\alpha+1}$ among the maps $v \in W^{1,p}_{\alpha+1}(\mathbf{B}^n,\mathbb{S}^{n-1})$ satisfying $v(y)=y$ on $\mathbb{S}^{n-1}$. 
\end{theorem}
The proof of this theorem relies on a construction which associates to each map $u:\mathbf{B}^n \rightarrow \mathbb{S}^{n-1}$ such
that $u(y)=y$ on $\mathbb{S}^{n-1}$, 
a map $\overline{u}~:\mathbf{B}^{n+1} \rightarrow \mathbb{S}^n$
such that $\overline{u}(x)=u(x)$ in $\mathbf{B}^n \times \{0\}$
and $\overline{u}(x)=x$ on the unit sphere $\mathbb{S}^n$, in such a way that, if
$u_0$ is the map defined in $\mathbf{B}^n$ by
$u_0(y)=\frac{y}{\|y\|}$, then $\overline{u}_0$ is exactly the map
defined on $\mathbf{B}^{n+1}$ by $\overline{u}_0(x)=\frac{x}{\|x\|}$. The energy  $E^{n+1}_{p,\alpha}(\overline{u})$ 
of  $\overline{u}$ is estimated above in terms of the energy $E^{n}_{p,\alpha+1}(u)$ of $u$ and the equality holds in the estimate for the map $u_0=\frac{y}{\|y\|}$. 

Thanks to Theorem \ref{th1} and the minimality results mentioned above, we deduce the following :

\begin{corollary}
The map $y/\|y\|$ minimizes ${E}^n_{p,\alpha}$ among the maps in $W_{\alpha}^{1,p}(\mathbf{B}^{n},\mathbb{S}^{n-1})$ which coincide with $y/\|y\|$ on $\mathbb{S}^{n-1}$ in the following cases~:

i) $\alpha \in \mathbb{N}$ and $p\in (n+\alpha-1,n+\alpha)$,

ii) $\alpha \geq \beta$ and $p$ is an integer in $\{1,\cdots, n+\beta-1\}$, for any  real number $\beta \geq 0$.

iii) $\alpha \in \mathbb{N}$, $n+\alpha \geq 7$ and $p\leq n+\alpha - 2\sqrt{n+\alpha -1}$,
%
%
%

%
\end{corollary}
%
%
%
\subsection{Construction of
$\overline{u}$ and proof of Theorem 1.1.}
%
%


Let $\mathbf {B}^{n+1}$ and $\mathbb{S}^n$ be the unit open ball and the
unit sphere of $\mathbb{R}^{n+1}$ and let $\mathbf{B}^{n}$ and
$\mathbb{S}^{n-1}$ be the unit open ball and the unit sphere of
$\mathbb{R}^n$ that we  identify with the subspace $\mathbb{R}^n= \mathbb{R}^n
\times\{0\}$ in $\mathbb{R}^{n+1}$. Moreover,  we write
$x=(x_1, \cdots, x_n, x_{n+1})$ a vector of $\mathbb{R}^{n+1}$,
$y=(y_1, \cdots ,y_n)$ a vector of $\mathbb{R}^n$,  $\langle .,. \rangle$ the
standard metric of $\mathbb{R}^{n+1}$ and $(e_1, \cdots, e_n,\,
e_{n+1})$ the standard basis of $\mathbb{R}^{n+1}$.
Let $\Pi_n$  be the projection defined by :
\begin{eqnarray*}
\Pi_n: \mathbf{B}^{n+1} & \longrightarrow& \mathbf{B}^n\\
 (x_1,\cdots,
x_{n+1})
                  &\longrightarrow& (x_1, \cdots, x_n, 0)=(x_1,\cdots, x_n).
\end{eqnarray*}
Consider the map $\varphi_n$ defined on $B^n \backslash
\mathbb{R}e_{n+1}$ by
$\varphi_n(x)=\frac{\Pi_n(x)}{\|\Pi_n(x)\|}$.\\
We define the map $\overline{u}$ defined by
\begin{eqnarray*}
\overline{u}: \mathbf{B}^{n+1} & \longrightarrow& \mathbf{B}^n\\ x
&\longrightarrow&  \langle e_{n+1},\frac{x}{\|x\|} \rangle e_{n+1}+
                 \langle \varphi_n(x),\frac{x}{\|x\|}\rangle u(\|x\|\varphi_n(x)).
\end{eqnarray*}
\begin{lem}
For any $x \in \mathbf{B}^{n+1}$, $\|\nabla
\overline{u}(x)\|^2=\frac{1}{\|x\|^2}+\|\nabla
u(\|x\|\varphi_n(x))\|^2$.

\end{lem}
\textbf{Proof.} 
For any $i\in \{1, \cdots, n\}$, we have, 
\begin{eqnarray*}
d\overline{u}(x).e_i &=& 
%
%
%
%
%
%
%
%
%
%
%
%
%
%
%
%
%
%
%
%
%
%
%
\frac{\langle \Pi_n(x),e_i\rangle}{\|x\|}\left(\frac{1}{\|\Pi_n(x)\|}-\frac{\|\Pi_n(x)\|}{\|x\|^2}\right)
u(\|x\|\varphi_n(x))\\
&&
-\frac{\langle e_{n+1},x\rangle \langle x,e_i\rangle }{\|x\|^3}e_{n+1}\\
&&
+\,du(\|x\|\varphi_n(x)).\Bigg(\frac{\langle \Pi_n(x),e_i\rangle \Pi_n(x)}{\|x\|^2}+e_i-
\frac{\langle \Pi_n(x),e_i\rangle \Pi_n(x)}{\|\Pi_n(x)\|^2}\Bigg) \\
\text{and,} &&\\
%
%
%
%
%
d\overline{u}(x).e_{n+1} &=&
%
%
%
%
%
%
%
\frac{e_{n+1}}{\|x\|}-\frac{\langle e_{n+1},x\rangle ^2}{\|x\|^3}e_{n+1}
-\frac{\langle e_{n+1},x \rangle }{\|x\|^3}\|\Pi_n(x)\|u(\|x\|\varphi_n(x))\\
&&
+\,\frac{\langle e_{n+1},x \rangle } {\|x\|^2}du(\|x\|\varphi_n(x)).\Pi_n(x).\\
\end{eqnarray*}
Since $\|u(x)\|^2=1$, one has
$\langle du(x).h,u(x)\rangle =0$ for any $h \in \mathbb{R}^n$. Hence, for
any $i \in \{1, \cdots,n\}$, we have,
\begin{eqnarray*}
\|d\overline{u}(x).e_i\|^2 &=&
\frac{\langle \Pi_n(x),e_i\rangle ^2}{\|x\|^2}\left(\frac{1}{\|\Pi_n(x)\|}
-\frac{\|\Pi_n(x)\|}{\|x\|^2}\right)^2
+\,\frac{\langle e_{n+1},x\rangle ^2\langle x,e_i\rangle ^2}{\|x\|^6}\\
&+&
\Big\|du(\|x\|\varphi_n(x)).\bigg(\langle \Pi_n(x),e_i\rangle \big(\frac{1}{\|x\|^2}
%
%
%
-\frac{1}{\|\Pi_n(x)\|^2}\big)\Pi_n(x)
+e_i\bigg)\Big\|^2.\\
\end{eqnarray*}
\begin{eqnarray*}
\|d\overline{u}(x).e_{n+1}\|^2 &=&
\frac{1}{\|x\|^2}\left(1-\frac{\langle e_{n+1},x\rangle ^2}{\|x\|^2}\right)^2+
\frac{\langle e_{n+1},x\rangle^2}{\|x\|^6}\|\Pi_n(x)\|^2\\
&&
+\,\frac{\langle e_{n+1},x\rangle ^2}{\|x\|^4}\|du(\|x\|\varphi_n(x)).\Pi_n(x)\|^2.
\end{eqnarray*}
Finally,  we have, 
\[
\|\nabla \overline{u}(x)\|^2=\frac{1}{\|x\|^2}+\|\nabla
u(\|x\|\varphi_n(x))\|^2. \quad\quad \blacksquare
\]

Let $x_{n+1}$  be a real number in $(0,1)$, consider the set
$A_{x_{n+1}}=(x_{n+1}e_{n+1}+e_{n+1}^{\bot}) \cap
\mathbf{B}^{n+1}\!\setminus\!\mathbb{R}e_{n+1}$, where $e_{n+1}^{\bot}$ is the orthogonal subspace to $\mathbb{R}e_{n+1}$ for $\langle .,. \rangle$. Let $\theta$ be
the map : 
\\
\begin{eqnarray*}
\theta: A_{x_{n+1}} & \longrightarrow & C_{x_{n+1}}=\{ y \in
\mathbb{R}^n; \|y\|>|x_{n+1}|\}\\ x=(x_1, \cdots, x_{n+1})
&\longrightarrow &  \|x\|\varphi_n(x)=y=(y_1, \cdots, y_n).\\
\end{eqnarray*}
\begin{lem}
For any $y \in \mathbb{R}^n$, the Jacobian determinant of $\theta^{-1}$ is : 
\[
Jac(\theta^{-1})(y)=
\frac{\big(\|y\|^2-x_{n+1}^2\big)^{\frac{n-2}{2}}}{\|y\|^{n-2}}.
\]
\end{lem}
\begin{dem}
For any $i \in \{1, \cdots,n\}$ and for any $x \in A_{x_{n+1}}$,
\begin{eqnarray*}
d\theta(x).e_i
&=&
%
%
%
%
%
%
\left(\frac{1}{\|x\|}-\frac{\|x\|}{\|\Pi_n(x)\|^2}
\right)\frac{\langle \Pi_n(x),e_i\rangle}{\|\Pi_n(x)\|}\Pi_n(x)+\frac{\|x\|}{\|\Pi_n(x)\|}e_i.\\
\end{eqnarray*}
Let us set, for any $i \in \{1, \cdots,n\}$,
\[
\lambda_i=\left(\frac{1}{\|x\|}-\frac{\|x\|}{\|\Pi_n(x)\|^2}
\right)\frac{\langle \Pi_n(x),e_i\rangle }{\|\Pi_n\|} \,\, \text{and}\,\,  \alpha=\frac{\|x\|}{\|\Pi_n(x)\|}.
\]
%
%
%
%
%
%
Hence, for any $i \in \{1, \cdots,n\}$,
\[
d \theta(x).e_i= \lambda_i\Pi_n(x)+\alpha e_i.
\]
Then, for any $x \in A_{x_{n+1}}$,
\begin{eqnarray*}
Jac(\theta)(x)&=&\det(\lambda_1\Pi_n(x)+\alpha
e_1,\lambda_2\Pi_n(x)+\alpha e_2, \cdots, \lambda_n\Pi_n(x)+\alpha
e_n)\\
&=&
\dmath\sum_{i=1}^n \det(\alpha e_1,\cdots,
\alpha e_{i-1},\lambda_i\Pi_n(x),\alpha e_{i+1}, \cdots, \alpha
e_n) +\alpha^n \det(e_1,\cdots,e_n)\\
%
%
%
%
%
&=&
\dmath\sum_{i=1}^n
\alpha^{n-1}\lambda_i\det(e_1,\cdots,e_{i-1},\Pi_n(x),e_{i+1},\cdots,e_n)
+ \alpha^n,\\
\end{eqnarray*}
where $\det$ is the determinant in the basis $(e_1, \cdots, e_{n
+1})$.
\begin{eqnarray*}
Jac( \theta)(x)&=&
\dmath\sum_{i=1}^n
\alpha^{n-1}\left(\frac{1}{\|x\|}-\frac{\|x\|}{\|\Pi_n(x)\|^2}
\right)\frac{\langle \Pi_n(x),e_i \rangle^2}{\|\Pi_n(x)\|}+ \alpha^n\\
&=&
\alpha^{n-1}\left(\frac{1}{\|x\|}-\frac{\|x\|}{\|\Pi_n(x)\|^2}
\right)\frac{\|\Pi_n(x)\|^2}{\|\Pi_n(x)\|}+ \alpha^n= \alpha^{n-2}.\\
%
%
%
%
%
%
%
%
\end{eqnarray*}
Therefore, we have
\[
Jac(\theta)(x)=\frac{\|x\|^{n-2}}{\|\Pi_n(x)\|^{n-2}}=\dmath
\frac{\|y\|^{n-2}}{\left(\|y\|^2-x_{n+1}^2\right)^{\frac{n-2}{2}}}.
\]
We deduce that,
\[
Jac(\theta^{-1})(y)=\dmath
\frac{\left(\|y\|^2-x_{n+1}^2\right)^{\frac{n-2}{2}}}{\|y\|^{n-2}}.
\]
\end{dem}
\begin{lem}
For any $x \in \mathbf{B}^{n+1}$, we have,
\begin{eqnarray*}
\int_{\mathbf{B}^{n+1}}\|x\|^{\alpha}\|\nabla \overline{u}(x)\|^{p}dx_1\cdots
dx_{n+1} &\leq&
n^{p/2-1}\int_{\mathbf{B}^{n+1}}\frac{1}{\|x\|^{p-\alpha}}dx_1\cdots
dx_{n+1}\\ &&+ 2 (1-1/n)^{1-p/2}\, W_{n-1} \int_{\mathbf{B}^n}
\|y\|^{\alpha+1}\|\nabla u(y)\|^p dy,\\
\end{eqnarray*}
where $W_{n-1}=\int_0^{\frac{\pi}{2}}(cos(\gamma))^{n-1}d \gamma.$
\end{lem}
%
%
%
%
\textbf{Proof.} Writing  $\frac{1}{\|x\|^2}+\|\nabla u(\|x\|\varphi_n(x))\|^2=\frac{1}{n}(n\frac{1}{\|x\|^2})+
(1-\frac{1}{n})(\frac{1}{1-\frac{1}{n}} \|\nabla u(\|x\|\varphi_n(x))\|^2)$ and using that $x \rightarrow x^{\frac{p}{2}}$ is a convex map,
for any $p \in [1,n+\alpha+1)$, we have,
\begin{eqnarray*}
\|\nabla \overline{u}(x)\|^{p} &=&
\left(\frac{1}{\|x\|^2}+\|\nabla
u(\|x\|\varphi_n(x))\|^2\right)^{p/2}\\
&\leq&
n^{p/2-1}\frac{1}{\|x\|^p}+(1-1/n)^{1-p/2}\|\nabla
u(\|x\|\varphi_n(x))\|^p.\\
\end{eqnarray*}
%
%
%
%
Hence,
\begin{eqnarray*}
\!\!\!\!\int_{\mathbf{B}^{n+1}}\|x\|^{\alpha}\|\nabla
\overline{u}(x)\|^{p}dx_1\cdots dx_{n+1} & \leq&
n^{p/2-1}\int_{\mathbf{B}^{n+1}}\|x\|^{\alpha}\frac{1}{\|x\|^p}dx_1\cdots
dx_{n+1}\\
%
%
%
+ (1-1/n\!\!\!\!&)^{1-p/2}&
\!\!\!\!\!\!\!\int_{\mathbf{B}^{n+1}}\|x\|^{\alpha}\|\nabla
u(\|x\|\varphi_n(x))\|^p dx_1\cdots dx_{n+1}\\
&\leq&
n^{p/2-1}\int_{\mathbf{B}^{n+1}}\|x\|^{\alpha}\frac{1}{\|x\|^p}dx_1\cdots
dx_{n+1}\\
+
\,2 (1-1/n)^{1-p/2}
\int_0^1\!\!\!&dx_{n+1}&\!\!\!\int_{\mathbf{B}^n}\|x\|^{\alpha}\|\nabla
u(\|x\|\varphi_n(x))\|^p dx_1\cdots dx_n.\\
\end{eqnarray*}
Using the change of variables $y=\theta(x)$ and Lemma 1.2 we
get, 
\begin{eqnarray*}
\int_{\mathbf{B}^{n+1}}\|x\|^{\alpha}\|\nabla \overline{u}(x)\|^{p}dx_1\cdots
dx_{n+1} &\leq&
n^{p/2-1}\int_{\mathbf{B}^{n+1}}\frac{1}{\|x\|^{p-\alpha}}dx_1\cdots
dx_{n+1}\\
+\,2 (1\!-\!1/n\!\!\!\!&)^{1-p/2}&\!\!\!\!
\int_0^1
\!\!\!dx_{n+1}\!\!\int_{C^{n+1}}\!\!\!\!\!\!\!\!\!
\frac{\left(\|y\|^2\!-\!x_{n+1}^2\right)^{\frac{n-2}{2}}}{\|y\|^{n-2}}\!\|y\|^{\alpha}\|\nabla
u(y)\|^{p}dy_1\!\cdots\! dy_n, 
\end{eqnarray*}
\begin{eqnarray*}
\!\!\!\!\int_{\mathbf{B}^{n+1}}\|x\|^{\alpha}\|\nabla
\overline{u}(x)\|^{p}dx_1\cdots dx_{n+1} &\leq&
n^{p/2-1}\int_{\mathbf{B}^{n+1}}\frac{1}{\|x\|^{p-\alpha}}dx_1\cdots
dx_{n+1}\\
+\,2 (1-\!1/n\!\!\!\!&)^{1-p/2}&\!\!\!\!
\int_{\mathbf{B}^n}\!\|y\|^{\alpha}\!\|\nabla u(y)\|^p
\!\!\left(\int_0^{\|y\|}
\frac{\left(\|y\|^2\!-\!x_{n+1}^2\right)^{\frac{n-2}{2}}}{\|y\|^{n-2}}
dx_{n+1}\!\!\right) \!\!dy.\\
\end{eqnarray*}
But we have,
\begin{eqnarray*}
\int_0^{\|y\|}
\left(\frac{\|y\|^2\!-\!x_{n+1}^2}{\|y\|^2}\right)^{\frac{n-2}{2}}
dx_{n+1} &=& \int_0^{\|y\|} \left(1-
\left(\frac{\!x_{n+1}}{\|y\|}\right)^2\right)^{\frac{n-2}{2}}dx_{n+1}\\
&=&
\int_0^1 \|y\|(1-t^2)^{\frac{n-2}{2}}dt =  \|y\|\int_0^{\pi
/2}\!\!\!\!(cos\gamma)^{n-1}\!d\gamma. \\
\end{eqnarray*}
Let us set $ W_{n-1}=\int_0^{\pi /2}(cos\gamma)^{n-1}d\gamma.$
Then, we have the inequality,
\begin{eqnarray*}
\int_{\mathbf{B}^{n+1}}\|x\|^{\alpha}\|\nabla \overline{u}(x)\|^{p}dx_1\cdots
dx_{n+1} &\leq&
n^{p/2-1}\int_{\mathbf{B}^{n+1}}\frac{1}{\|x\|^{p-\alpha}}dx_1\cdots
dx_{n+1}\\
&&
+ \,2 (1-1/n)^{1-p/2}\, W_{n-1} \int_{\mathbf{B}^n} \|y\|^{\alpha+1}\|\nabla
u(y)\|^p dy. \quad \blacksquare\\ 
\end{eqnarray*}
\begin{lem}
Let $\Gamma$ be the Gamma function defined by,
\[
\Gamma(x) = \int_0^{+ \infty} t^{x-1} e^{-t} dt.
\]
We have,
\[
W_{n-1} \frac{\Gamma(\frac{n+1}{2})}{\Gamma(\frac{n}{2})}=
\frac{\sqrt{\pi}}{2}.
\]
\end{lem}
\begin{dem}
From the equality $\Gamma(x+1)=x\Gamma(x)$ for any $x \in
(0,+\infty)$, we have, if $n=2l$, where $l \in \mathbb{N}$,
%
%
%
%
%
\[
\Gamma(\frac{2l+1}{2})= \frac{2l-1}{2}\frac{2l-3}{2}\cdots
\frac{3}{2}\frac{1}{2}\Gamma(\frac{1}{2})
%
%
\]
and
%
%
%
%
%
\[
\Gamma(\frac{2l}{2})= \frac{2l-2}{2}\frac{2l-4}{2}\cdots
\frac{2}{2}\Gamma(1).
\]
Moreover we have,
%
\[
W_{2l-1}=\frac{(2l-2)(2l-4)\cdots 2}{(2l-1)(2l-3)\cdots 3.1}
\]
%
%
then,
\[
W_{2l-1} \frac{\Gamma(\frac{2l+1}{2})}{\Gamma(\frac{2l}{2})}=
\frac{\sqrt{\pi}}{2}.
\]
If $n=2l+1$ where $l \in \mathbb{N}$,
we obtain,
\[
\Gamma(\frac{2l+2}{2})=
\frac{2l}{2}\frac{2l-2}{2}\cdots\frac{2}{2}\Gamma(1),
\]
\[
\Gamma(\frac{2l+1}{2})= \frac{2l-1}{2}\cdots\frac{3}{2}
\frac{1}{2}\Gamma(\frac{1}{2})
\]
and
\[
W_{2l}=\frac{(2l-1)(2l-3)\cdots 3.1}{(2l)(2l-2)\cdots
2}.\frac{\pi}{2}.
\]
We deduce that
\[
W_{2l} \frac{\Gamma(\frac{2l+2}{2})}{\Gamma(\frac{2l+1}{2})}=
\frac{\sqrt{\pi}}{2}. \quad \quad \blacksquare 
\]
\end{dem}
\begin{dem}
\textbf{of  Theorem 1.1.}
Suppose that $\overline{u}_0$ is a minimizer map of $E_{p,\alpha}$. 
Since the measure of $\mathbb{S}^{n}$ is $|\mathbb{S}^{n}|=\frac{2\pi^{\frac{n+1}{2}}}{\Gamma(\frac{n+1}{2})}$ and since $\|\nabla \overline{u}_0\|(x)=\frac{n^{\frac{p}{2}}}{\|x\|}$,
 by Lemma 1.3, for any map $u \in W_{\alpha+1}^{1,p}(\mathbf{B}^n,
\mathbb{S}^{n-1})$ satisfying $u(y)=y$ on $\mathbb{S}^{n-1}$, we
have,
\begin{eqnarray*}
\frac{n^{p/2}}{n+1+\alpha-p}\frac{2
\pi^{\frac{n+1}{2}}}{\Gamma(\frac{n+1} {2})}
&=&
\int_{\mathbf{B}^{n+1}}\|x\|^{\alpha}\|\nabla \overline{u_0}(x)\|^{p}dx_1\cdots
dx_{n+1}\\
 &\leq& \int_{\mathbf{B}^{n+1}}\|x\|^{\alpha}\|\nabla
\overline{u}(x)\|^{p}dx_1\cdots dx_{n+1}\\
&\leq&
n^{p/2-1}\int_{\mathbf{B}^{n+1}} \frac{1}{\|x\|^{p-\alpha}}dx_1 \cdots
dx_{n+1}\\
&&
\hspace{1cm}+ 2 (1\!\!-\!\!1/n)^{1-p/2}\, W_{n-1}
\!\!\int_{\mathbf{B}^n} \|y\|^{\alpha+1}\nabla u(y)\|^p dy\\
&\leq &
\frac{1}{n} \frac{n^{p/2}}{n\!+\!1\!+\alpha-\!p}
\frac{2\pi^{\frac{n+1}{2}}}{\Gamma(\frac{n+1}{2})}\\
&&
\hspace{1cm} +\, 2W_{n-1}
(1-\frac{1}{n})(\frac{n}{n-1})^{p/2}\int_{\mathbf{B}^n}
\|y\|^{\alpha+1}\|\nabla u(y)\|^p dy.\\
\end{eqnarray*}
Then we have,
\begin{eqnarray*}
2W_{n-1}
\left(1\!-\!\frac{1}{n}\right)\left(\frac{n}{n\!-\!1}\right)^{p/2}
\dmath\frac{\int_{\mathbf{B}^n} \|y\|^{\alpha+1}\|\nabla u(y)\|^p
dy}{\int_{\mathbf{B}^n} \|y\|^{\alpha+1}\|\nabla \frac{y}{\|y)\|}\|^p dy}
&\geq&
%
%
%
%
%
%
%
\frac{n^{p/2}}{n\!+\!1\!+ \alpha-\!p}\,\frac{2
\pi^{\frac{n+1}{2}}}{\Gamma(\frac{n+1} {2})}\\
&&
\times\frac{n\!+\!1\!+\alpha-\!p}{(n-1)^{p/2}}\frac{\Gamma(\frac{n}
{2})}{2 \pi^{\frac{n}{2}}}\left(1\!-\!\frac{1}{n}\right)\\
\end{eqnarray*}
%
%
%
%
%
%
%
%
%
%
By Lemma 1.4 we finally get,
\[
\int_{\mathbf{B}^n} \|y\|^{\alpha+1}\|\nabla u(y)\|^p dy \geq
\int_{\mathbf{B}^n} \|y\|^{\alpha+1}\|\nabla \frac{y}{\|y\|}\|^p dy \quad
\quad \quad \blacksquare
\]
\end{dem}
\end{document}